\documentclass[12]{article}
\usepackage{amsfonts}
\usepackage{theorem}
\newtheorem{definition}{Definition}
\newtheorem{theorem}{Theorem}
\newtheorem{corollary}{Corollary}
\newtheorem{lemma}{Lemma}
\newtheorem{remark}{Remark}
\usepackage[unicode,colorlinks,plainpages=false]{hyperref}
\newcommand{\openbox}{$\begin{array}{c}
\hspace*{-0.55em}\sqcap \hspace*{-0.60em}\\[-0.4em] \hline
\multicolumn{1}{c}{\hspace*{-0.60em}}\\[-0.8em]
\end{array}
$}

\usepackage{enumerate}

\begin{document}

\centerline{\bf The separator of a subset of a semigroup}

\bigskip

\bigskip

\centerline{\bf Attila Nagy}

 \medskip
\centerline{Department of Algebra}
\centerline{Mathematical Institute}
\centerline{Budapest University of Technology and Economics}
\centerline{e-mail: nagyat@math.bme.hu}

\bigskip
\small

\centerline{\bf 1. Introduction}

\bigskip

In this paper we introduce a new notion by the help of the idealizer. This new notion is the separator of a subset of a semigroup. We investigate the properties of the separator in an arbitrary semigroup and characterize the unitary subsemigroups and the prime ideals by the help of their separator. We give conditions which imply that a maximal ideal is prime. The last section of this paper treats the separator of a free subsemigroup of a free semigroup.

\medskip

If $A$ is a subset of a semigroup $S$, then $\overline{A}$ denotes the subset $S\setminus A$. For other symbols we refer to \cite{2}.

\bigskip

\bigskip

\centerline{\bf 2. The definition and basic properties of the separator}

\bigskip

Let $S$ be a semigroup and $A$ any subset of $S$. As known \cite{1}, the idealizer of $A$ is the set of all elements $x$ of $S$ which satisfy the following conditions: $Ax\subseteq A$, $xA\subseteq A$. The idealizer of a subset $A$ is denoted by $IdA$.

If $A$ is an empty set, then $IdA$ equals to $S$. It is evident that $IdA$ is either empty or a subsemigroup of $S$.

\begin{definition}\label{df1} Let $S$ be a semigroup and $A$ any subset of $S$. Then $IdA\cap Id\overline{A}$ will be called the separator of $A$ and denoted by $SepA$.
\end{definition}

\medskip

In other words: An element $x$ of $S$ belongs to $SepA$ ($A\subseteq S$) if and only if $xA\subseteq A$, $Ax\subseteq A$, $x\overline{A}\subseteq \overline{A}$, $\overline{A}x\subseteq \overline{A}$.

\begin{remark}\label{rm1} For any subset $A$ of $S$, $SepA$ is either empty or a subsemigroup of $S$. Moreover, $SepA=Sep\overline{A}$. In particular, $Sep\emptyset =SepS=S$.
\end{remark}

\begin{remark}\label{rm2} If $S$ is a semigroup with an identity element, then the identity element belongs to $SepA$ for any subset $A$ of $S$.
\end{remark}

\begin{remark}\label{rm3} If $R$ is an ideal of a semigroup $S$ with $R\neq S$, then $R\cap SepR=\emptyset$.
\end{remark}

\begin{theorem}\label{th1} Let $\{A_f:\ f\in F\}$ be any non-empty family of subsets of a semigroup $S$. Then
$\cap _{f\in F}SepA_f\subseteq Sep(\cup _{f\in F}A_f)$ and
$\cap _{f\in F}SepA_f\subseteq Sep(\cap _{f\in F}A_f)$.
\end{theorem}

\newpage

{\bf Proof}. Let $t\in \cap _{f\in F}SepA_f$. Then $t[\cup _{f\in F}A_f]=\cup _{f\in F}tA_f\subseteq \cup _{f\in F}A_f$ and
$t\overline{\cup _{f\in F}A_f}=t[\cap _{f\in F}\overline{A}_f]\subseteq \cap _{f\in F}t\overline{A}_f\subseteq \cap _{f\in F}\overline{A}_f=\overline{\cup _{f\in F}A_f}$.

Similarly, $[\cup _{f\in F}A_f]t\subseteq \cup _{f\in F}A_f\quad \hbox{and}\quad \overline{\cup _{f\in F}A_f}t\subseteq \overline{\cup _{f\in F}A_f}$. Consequently,\newline $t\in Sep(\cup _{f\in F}A_f)$ and the first part of the theorem is proved.

As $\cap _{f\in F}SepA_f=\cap _{f\in F}Sep \overline{A}_f\subseteq Sep(\cup _{f\in F}\overline{A}_f)=Sep (\overline{\cap _{f\in F}A_f})=Sep(\cap _{f\in F}A_f)$, the theorem is proved. \hfill\openbox

\begin{corollary}\label{cr1} $SepA\cap SepSepA\subseteq Sep(A\cup SepA)$ for any subset $A$ of a semigroup.
\end{corollary}

\begin{theorem}\label{th2} If $A$ is a subsemigroup of a semigroup, then $A\cup SepA$ is so.
\end{theorem}

{\bf Proof}. Let $A$ be a subsemigroup of a semigroup $S$. We may assume that $SepA\neq \emptyset$. Let $x, y\in A\cup SepA$. Since $SepA$ is also a subsemigroup of $S$, we have to consider only the case when $x\in A$ and $y\in SepA$. Then $xy$ and $yx$ belong to $A$ by the definition of the separator.\hfill\openbox

\begin{theorem}\label{th3} If $A$ is a subset of a semigroup $S$ such that $SepA\neq \emptyset$, then either $SepA\subseteq A$ or $SepA\subseteq \overline{A}$.
\end{theorem}

{\bf Proof}. Let $A$ be a subset of a semigroup $S$ and assume that $A\cap SepA\neq \emptyset$. Let $x$ be an arbitrary element of $SepA$ and $a$ an element of $A\cap SepA$. Then $xa\in A$. If $x$ were in $\overline{A}$, then $xa$ would be in $\overline{A}$, too (because $a\in sepA$), and this would be a contradiction. Consequently, if $A\cap SepA\neq \emptyset$ and $x\in SepA$, then $x\in A$. In other words, $A\cap SepA\neq \emptyset$ implies $SepA\subseteq A$.\hfill\openbox

\begin{theorem}\label{th4} Let $\phi$ be a homomorphism of a semigroup $S$ onto itself and $R_1$, $R_2$ subsemigroups of $S$ such that $\phi ^{-1}(R_2)=R_1$. Then $\phi (SepR_1)=SepR_2$.
\end{theorem}

{\bf Proof}. Let $x$ be an arbitrary element of $SepR_1$. We prove that $\phi (x)$ belongs to $SepR_2$. For this purpose, let $y$ be an arbitrary element of $R_2$. Then there exists a $y_0\in S$ such that $\phi (y_0)=y$. Assume $\phi (x)y\in \overline{R}_2$. Then
$\phi (xy_0)=\phi (x)\phi (y_0)\in \overline{R}_2$ and therefore $xy_0\in \overline{R}_1$. Thus $y_0\in \overline{R}_1$, because $x\in SepR_1$. Hence $y=\phi (y_0)\in \overline{R}_2$ which contradicts the assumption that $y\in R_2$. Consequently, $\phi (x)y$ and, similarly, $y\phi (x)$ are in $R_2$ for every $x\in SepR_1$, $y\in R_2$.

Now let $z$ be an arbitrary element of $\overline{R}_2$. Then there exists a $z_0\in S$ such that $\phi (z_0)=z$. Assume $\phi (x)z\in R_2$. Then $\phi (xz_0)=\phi (x)\phi (z_0)$ belongs to $R_2$ and therefore $xz_0\in R_1$. Thus $z_0\in R_1$, because $x\in SepR_1$. Hence $z=\phi (z_0)\in R_2$ which contradicts the assumption that $z\in \overline {R}_2$. Consequently $\phi (x)z\in \overline{R}_2$ and, similarly, $z\phi (x)\in \overline{R}_2$. Thus $\phi (x)$ belongs to $SepR_2$, indeed. Hence $\phi (SepR_1)\subseteq SepR_2$.

We complete the proof by showing that $t\notin SepR_1$ implies $\phi (t)\notin SepR_2$. In fact, if $t\notin SepR_1$, then either there exists an element $u$ in $R_1$ such that at least one of the inclusions $tu\in \overline{R}_1$ and $ut\in \overline{R_1}$ holds or there exists an element $v\in \overline{R}_1$ such that at least one of $tv\in R_1$ and $vt\in R_1$ holds. Consider the case when $u\in R_1$ and $tu\in \overline{R}_1$ (the other cases can be discussed similarly). Then $\phi (u)\in R_2$ and $\phi (t)\phi (u)\in \overline{R}_2$ which means that $\phi (t)\notin SepR_2$, indeed.\hfill\openbox

\begin{corollary}\label{cr2} If $\phi$ is an isomorphism of a semigroup $S$ onto itself and $A$ and $B$ are subsemigroups of $S$ such that $\phi (A)=B$, then $\phi (SepA)=SepB$.
\end{corollary}

\newpage

\centerline{\bf 3. Separator including and separator excluding subsets}

\medskip

\begin{definition}\label{df2} A subset $A$ of a semigroup is said to be separator including [excluding] if $SepA\subseteq A$ [$Sep A \subseteq \overline{A}$]. In the case of $SepA=\emptyset$, the subset $A$ will be considered separator including as well as separator excluding.
\end{definition}

\begin{remark}\label{rm4} If a subset $A$ of a semigroup is separator including [excluding], then, evidently, $\overline{A}$ is separator excluding [including]. In particular, the semigroup $S$ is a separator including subset of itself and the empty set is a separator excluding one.
\end{remark}

\begin{remark}\label{rm5} If $A$ and $B$ are separator including subsemigroups of a semigroup $S$ and $\phi$ is a homomorphism of $A$ onto $B$, then, in general, $\phi (SepA)\neq SepB$. But if $\phi$ is a homomorphism of $A$ onto $B$ such that $SepA=\phi ^{-1}(SepB)$, then $\phi$ is a homomorphism of $SepSepA$ onto $SepSepB$.
\end{remark}

\medskip

{\bf Example}:
Consider the semigroup $S=\{ 0, a, b, 1\}$ in which the operation is given by the following Cayley table:

\begin{table}[ht]
 \center{\begin{tabular}{l|l l l l }
 &$0$&$a$&$b$&$1$\\ \hline
$0$&$0$&$0$&$0$&$0$\\
$a$&$0$&$a$&$0$&$a$\\
$b$&$0$&$0$&$b$&$b$\\
$1$&$0$&$a$&$b$&$1$\\
 \end{tabular}}
\end{table}

In this semigroup, $Sep\{ 1\}=\{ 1\}$, $Sep\{a, b, 1\}=\{ a, b, 1\}$. As $Sep\{ a\}=\{ 1\}$, \newline $\{ a\}\cap Sep\{ a\}=\emptyset$ and
$\{ a\} \cup Sep\{ a\}\neq S$. Since $Sep\{ 0, a, b\}=\{ 1\}$, we have \newline $\{ 0, a, b\}\cap Sep\{ 0, a, b\}=\emptyset$ and $\{ 0, a, b\}\cup Sep\{ 0, a, b\}=S$.

\begin{theorem}\label{th5} Let $A$ be a separator including [excluding] subset and $B$ an arbitrary subset of a semigroup $S$ such that
$SepA\cap SepB\neq \emptyset$. Then $A\cup B$ [$A\cap B$] also is a separator including [excluding] subset of $S$ and its separator is non-empty.
\end{theorem}

{\bf Proof}. By Theorem~\ref{th3}, our assertion concerning a separator including $A$ will be proved if we show that the intersection of $A\cup B$ and $Sep(A\cup B)$ is non-empty. Using Theorem~\ref{th1} and the condition $SepA \subseteq A$, we get $(A\cup B)\cap Sep(A\cup B)\supseteq A\cap SepA \cap SepB=SepA\cap SepB$ and the last term is non-empty by one of the conditions. Consequently $(A\cup B)\cap Sep(A\cup B)\neq \emptyset$, indeed.

Consider the assertion concerning a separator excluding $A$. Then $\overline{A}$ is separator including. Since $Sep\overline{A}\cap Sep\overline{B}=
SepA\cap SepB\neq \emptyset$, by our assertion proved just now, $\overline{A}\cup \overline{B}$ is separator including and $Sep(\overline{A}\cup \overline{B})\neq \emptyset$. Consequently, $A\cap B=\overline{\overline{A}\cup \overline{B}}$ is separator including and $\emptyset \neq Sep(\overline{A}\cup \overline{B})=Sep(\overline{\overline{A}\cup \overline{B}})=Sep(A\cap B)$.\hfill\openbox

\begin{corollary}\label{cr3} Let $S$ be a semigroup and $A, B\subseteq S$. If $A$ is a separator including subset and $B$ is a separator excluding one such that $SepA\cap SepB\neq \emptyset$, then $A\cup B$ is separator including, $A\cap B$ separator excluding and $Sep(A\cup B)\neq \emptyset$, $Sep(A\cap B)\neq \emptyset$.
\end{corollary}

\newpage

\begin{theorem}\label{th6} Let $A$ and $B$ separator including [excluding] subsets of a semigroup $S$ such that $SepA\cap SepB\neq \emptyset$. Then $A\cap B$ [$A\cup B$] is also a separator including [excluding] subset of $S$ and its separator is non-empty.
\end{theorem}

{\bf Proof}. By Theorem~\ref{th3}, our assertion concerning a separator including $A$ and $B$ will be proved if we show that the intersection of $A\cap B$ and $Sep(A\cap B)$ is non-empty. Using Theorem~\ref{th1} and the condition that $SepA\subseteq A$ and $SepB\subseteq B$, we get
$A\cap B\cap Sep(A\cap B)=A\cap B\cap Sep(\overline{A\cap B})=A\cap B\cap Sep(\overline{A}\cup \overline{B})\supseteq A\cap B\cap Sep\overline{A}\cap Sep\overline{B}=
A\cap B\cap SepA\cap SepB=SepA\cap SepB$ and the last term is non-empty by one of the conditions. Consequently, $A\cap B\cap Sep(A\cap B)\neq \emptyset$, indeed.

Consider the assertion concerning separator excluding $A$ and $B$. Then $\overline{A}$ and $\overline{B}$ are separator including. Since $Sep\overline{A}\cap Sep\overline{B}=SepA\cap SepB\neq\emptyset$, $\overline {A}\cap\overline{B}$ is separator including and $Sep(\overline {A}\cap \overline{B})\neq \emptyset$. Consequently, $A\cup B=\overline{\overline{A}\cap \overline{B}}$ is separator excluding and
$\emptyset \neq Sep(\overline{A}\cap\overline{B})=Sep\overline{\overline{A}\cap \overline{B}}=Sep(A\cup B)$.\hfill\openbox

\begin{corollary}\label{cr4} Let $A$ be a separator including subset and $B$ an arbitrary subset of a semigroup $S$. If $SepA\cap SepB\neq \emptyset$, then $B\setminus A$ is separator excluding and $Sep(B\setminus A)\neq \emptyset$, where $B\setminus A$ denotes $B\cap \overline{A}$. Obviously, $B=(B\setminus A)\cup (B\cap A)$ holds for every subsets $A$ and $B$ of $S$. Hence, if $A$ is separator including, $B$ is separator excluding and $SepA\cap SepB\neq \emptyset$, then $B$ is a union of two separator excluding subsets (because $B\setminus A=B\cap \overline{A}$ and $A\cap B$ are separator excluding).
\end{corollary}

\begin{theorem}\label{th7} The separator of any subset $A$ of a semigroup is separator including (that is, $SepSepA\subseteq SepA$), provided that
$SepA\neq \emptyset$.
\end{theorem}

{\bf Proof}. Let $A$ be a subset of a semigroup. We may assume that $SepSepA\neq \emptyset$. Let $t$ be an arbitrary element of $SepSepA$. We prove
$t\in SepA$. Let $a$ be an arbitrary element of $A$. Then, for every $x\in SepA$, $xt\in SepA$ and thus $xta\in A$. If $ta$ was in $\overline{A}$, then $xta$ would be in $\overline{A}$, because $x\in SepA$, and this would be a contradiction. Hence $ta\in A$. Similarly, $at\in A$. Now, let $b$ be an arbitrary element of $\overline{A}$. We prove that $tb\in \overline{A}$. Let $x$ be an arbitrary element of $SepA$. Then $xt\in SepA$ and so $xtb\in \overline{A}$. If $tb$ was in $A$, then $xtb$ would be in $A$, because $x\in SepA$, and this would be a contradiction. Hence $tb\in \overline{A}$, indeed. Similarly, $bt\in \overline{A}$. Thus, by the definition of the separator, $t$ belongs to $SepA$. This implies that
$SepSepA\subseteq SepA$.\hfill\openbox

\begin{corollary}\label{cr5} If the separator of a subset $A$ of a semigroup $S$ is a minimal subsemigroup, then $SepSepA$ either is empty or equal to $SepA$.
\end{corollary}

\bigskip

\centerline{\bf 4. Unitary subsemigroups and prime ideals}

\medskip

\begin{definition}\label{df3} A subsemigroup $U$ of a semigroup $S$ is called unitary if $ab, a\in U$ implies $b\in U$ and $ab, b\in U$ implies $a\in U$. In other words: A subsemigroup $U$ of $S$ is unitary if $ab\in U$ implies either $a, b\in U$ or $a, b\in \overline{U}$.
\end{definition}

\begin{theorem}\label{th8} For a subsemigroup $A$ of a semigroup $S$, the following two assertions are equivalent:
\begin{enumerate}
\item[(i)] $A=SepA$,
\item[(ii)] $A$ is a unitary subsemigroup of $S$.
\end{enumerate}
\end{theorem}

\newpage

{\bf Proof}. $(i)$ implies $(ii)$: Let $A$ be a subsemigroup of a semigroup $S$ such that $A=SepA$. If $a\in A\ (=SepA)$ and $b\in \overline{A}$ or $a\in \overline{A}$ and $b\in A$, then $ab\in \overline{A}$. Consequently, $ab\in A$ implies either $a, b\in A$ or $a, b\in \overline{A}$. Thus $A$ is a unitary subsemigroup of $S$.

$(ii)$ implies $(i)$: Let $A$ be a unitary subsemigroup of a semigroup $S$ and $a\in A$. Then $aA\subseteq A$ and $Aa\subseteq A$, because $A$ is a subsemigroup of $S$. Moreover $a\overline{A}\subseteq \overline{A}$ and $\overline{A}a\subseteq \overline{A}$, because $A$ is unitary in $S$.
Hence $A\subseteq SepA$. Consequently, $A=SepA$ by Theorem~\ref{th3}.\hfill\openbox

\medskip

The following two theorems deal with the prime ideals of semigroups. The first of them, which is partly known (\cite{2}), characterizes the prime ideals, and the second gives conditions for a maximal ideal to be prime.

\medskip

An ideal $P$ of a semigroup $S$ is called a prime ideal if $P\neq S$ and, for every $a, b\in S$, $ab\in P$ implies $a\in P$ or $b\in P$.

\begin{theorem}\label{th9} A subsemigroup $P$ of a semigroup $S$ is a prime ideal of $S$ if and only if the complement of $P$ is an unitary subsemigroup of $S$.
\end{theorem}

{\bf Proof}. For any ideal $P\neq S$ of a semigroup $S$, $SepP\subseteq \overline{P}$ by Remark~\ref{rm3}. Moreover, $\overline{P}P\subseteq P$ and $P\overline{P}\subseteq P$. Assume that the ideal $P$ is prime in $S$. Then also $\overline{P}\ \overline{P}\subseteq \overline{P}$ which implies
(together with $\overline{P}P\subseteq P$ and $P\overline{P}\subseteq P$) that
 $\overline{P}\subseteq SepP=Sep\overline{P}$. By Theorem~\ref{th3}, $\overline{P}=Sep\overline{P}$. By Theorem~\ref{th8}, $\overline{P}$ is a unitary subsemigroup of $S$.

Conversely, assume that $P$ is a subsemigroup of semigroup $S$ and $\overline{P}$ is a unitary subsemigroup of $S$. Then $\overline{P}=Sep\overline{P}=SepP$. Let $x$ be an arbitrary element of $S$. If $x\in \overline{P}$, then $px\in P$ and $xp\in P$ ($p\in P$) because $\overline{P}=SepP$. If $x\in P$ then $xP\subseteq P$ and $Px\subseteq P$, because $P$ is a subsemigroup. Thus $xP\subseteq P$ and $Px\subseteq P$. Consequently, $P$ is an ideal of $S$. Since $\overline{P}$ is a subsemigrop by the one of the conditions, the ideal $P$ is prime (\cite{2}).\hfill\openbox

\medskip

An ideal $M\neq S$ of a semigroup $S$ is called a maximal ideal if, for every ideal $A$ of $S$, the assumption $M\subseteq A\subseteq S$ implies $M=A$ or $A=S$.

\begin{theorem}\label{th10} Let $I$ be a maximal ideal of a semigroup $S$. If $I$ is a maximl subsemigroup of $S$ and $SepI\neq \emptyset$, then
$I$ is a prime ideal.
\end{theorem}

{\bf Proof}. Let $I$ be a maximal ideal of $S$. By Remark~\ref{rm3}, $I\cap SepI =\emptyset$ and so, by Theorem~\ref{th2}, $I\cup SepI$ is a subsemigroup of $S$ such that $I\cup SepI\supset I$. As $I$ is also a maximal subsemigroup of $S$, we get $I\cup SepI=S$ and so $\overline{I}=SepI=Sep\overline{I}$. By Theorem~\ref{th8}, $\overline{I}$ is a unitary subsemigroup of $S$. By Theorem~\ref{th9}, $I$ is a prime ideal.\hfill\openbox

\bigskip

\centerline{\bf 5. Separators in a free semigroup}

\bigskip

In this section we shall use the next lemma (see \cite{2}, \S 9.1) several times:

\begin{lemma}\label{lm1} A subsemigroup $T$ of a free semigroup $S$ is a free subsemigroup if and only if $sT\cap T\neq \emptyset$ and $Tt\cap T\neq \emptyset$ together imply $s\in T$ for each element $s$ of $S$.
\end{lemma}

\begin{theorem}\label{th11} Any free subsemigroup of a free semigroup is separator including.
\end{theorem}

{\bf Proof}. Let $T$ be a free subsemigroup of a free semigroup $S$. We may assume that $SepT\neq \emptyset$. Let $s$ be an arbitrary element of $SepT$. Then $sT\subseteq T$ and $Ts\subseteq T$. Thus $sT\cap T\neq \emptyset$, $Ts\cap T\neq \emptyset$ and so $s\in T$ by the Lemma~\ref{lm1}. Consequently,
$SepT\subseteq T$.\hfill\openbox

\newpage

\begin{theorem}\label{th12}  The separator of any free subsemigroup $T$ of a free semigroup is free, unless $SepT=\emptyset$.
\end{theorem}

{\bf Proof}. Assume hat $T$ is a free subsemigroup of a free semigroup $S$ and $SepT\neq \emptyset$. Let $s$ be an arbitrary element of $S$ such that $s(SepT)\cap SepT\neq \emptyset$ and $(SepT)s\cap SepT\neq \emptyset$. Then $sT\cap T\neq \emptyset$ and $Ts\cap T\neq \emptyset$, becuse $SepT\subseteq T$ by Theorem~\ref{th11}. Hence $s\in T$, because $T$ is a free subsemigroup of $S$. We prove that $s\in SepT$. First, $sT\subseteq T$ and $Ts\subseteq T$, because $s\in T$. Next, let $t$ be an arbitrary element of $\overline{T}$. We have to show that $st$ and $ts$ belong to $\overline{T}$. By the condition $(SepT)s\cap SepT\neq \emptyset$, there exist $m_1, m_2\in SepT\ (\subseteq T)$ such that $m_1s=m_2$. Hence $m_2t\in \overline{T}$. Now, if we suppose $st\in T$, then we get $m_2t=(m_1s)t=m_1(st)\in T$ which contradicts $m_2t\in \overline{T}$. Consequently, $s\overline{T}\subseteq \overline{T}$ and, similarly, $\overline{T}s\subseteq \overline{T}$. Thus $s\in SepT$, as it was asserted. Thus we have shown that the condition
$s(SepT)\cap SepT\neq \emptyset$ and $(SepT)s\cap SepT\neq \emptyset$ imply that $s\in SepT$ for each $s\in S$. It follows (see again Lemma~\ref{lm1}) that $SepT$ is a free subsemigroup of $S$.\hfill\openbox

\begin{theorem}\label{th13} For any subset $T\neq \emptyset$ of a free semigroup $S$, the following two conditions are equivalen.
\begin{enumerate}
\item[(i)] $T$ is a free subsemigroup with the property that, for any elements $t\in T$ and $s\in S$, the condition $tst\in T$ implies both
$ts\in T$ and $st\in T$.
\item[(ii)] $T=SepT$.
\end{enumerate}
\end{theorem}

{\bf Proof}. Assume that $(i)$ holds. By Theorem~\ref{th11}, it is sufficient to show that $T\subseteq SepT$. Let $t$ be an arbitrary element of $T$. Then $tT\subseteq T$ and $Tt\subseteq T$, because $T$ is a subsemigroup. Let $s$ be an arbitrary element of $\overline{T}$. We show that
$ts\in \overline{T}$ and $st\in \overline{T}$. If we assume that $ts\in T$ or $st\in T$, then $tst\in T$ and so, by $(i)$, $ts\in T$ and $st\in T$. Hence $Ts\cap T\neq \emptyset$ and $sT\cap T\neq \emptyset$. As $T$ is a free subsemigroup, we get $s\in T$ by Lemma~\ref{lm1}. This contradicts
$s\in \overline{T}$. Consequently, $t\overline{T}\subseteq \overline{T}$ and $\overline{T}t\subseteq \overline{T}$. Thus $t\in SepT$ and so $T\subseteq SepT$.

Conversely, let $T$ be a non-empty subset of a free semigroup $S$ such that $T=SepT$. Then $T$ is a subsemigroup of $S$ by Remark~\ref{rm1}. Let $s$ be an arbitrary element of $S$ such that $sT\cap T\neq \emptyset$ and $Ts\cap T\ne \emptyset$. Then $s\in T$, because if $s\in \overline{T}$ then $sT=s(SepT)=s(Sep\overline{T})\subseteq \overline{T}$ and so $sT\cap T=\emptyset$; this is a contradiction. Hence $T$ is a free subsemigroup by Lemma~\ref{lm1}. Let $t\in T$, $s\in S$ be an arbitrary elements such that $tst\in T$. Since $T=SepT$, then $T$ is a unitary subsemigroup of $S$ (see Theorem~\ref{th8}) and so $t\in T$ implies $st, ts\in T$.\hfill\openbox

 \begin{theorem}\label{th14} Let $S$ be a free semigroup and $A$ a subsemigroup of $S$ such that $SepA\supseteq A^n$ for some positive integer $n$. Then $A$ is a free subsemigroup.
 \end{theorem}

 {\bf Proof}. Let $s\in S$ be arbitrary with $sA\cap A\neq \emptyset$ and $As\cap A\neq \emptyset$. Then there is an element $a\in A$ such that $sa\in A$. Since $A$ is a subsemigroup, $sa^n\in A$. Hence $s\in A$, because $a^n\in A^n\subseteq SepA$. By Lemma~\ref{lm1}, $A$ is a free subsemigroup of $S$.\hfill\openbox

\end{document}